\documentclass[11pt]{amsart}
\usepackage{amsfonts}
\usepackage{latexsym,amsmath,amsthm,amssymb,amsfonts}

\numberwithin{equation}{section}
\allowdisplaybreaks

\def\endproof{$\hfill\Box$\\}
\def\s{\,\,\,\,}
\def\R{\mathbb{R}}
\def\C{\mathbb{C}}
\def\CI{W^{2,2}_{{\rm conf}}}
\def\N{\mathbb{N}}
\def\mv{2.0ex}
\def\Re{\mbox{Re}}
\numberwithin{equation}{section}
\newtheorem{theorem}{Theorem}[section]
\newtheorem{lem}[theorem]{Lemma}
\newtheorem{thm}[theorem]{Theorem}
\newtheorem{pro}[theorem]{Proposition}
\newtheorem{cor}[theorem]{Corollary}
\newtheorem{defi}[theorem]{Definition}
\newtheorem{rem}[theorem]{Remark}

\def\dint{\displaystyle{\int}}

\begin{document}
\title[Weak limit of an immersed surface sequence
with bounded Willmore functional]
{\bf Weak limit of an immersed surface sequence
with bounded Willmore functional}
\author[Yuxiang Li]{Yuxiang Li\\
{\small\it Department of Mathematical Sciences},\\
{\small\it Tsinghua University,}\\
{\small\it Beijing 100084, P.R.China.}\\
{\small\it Email: yxli@math.tsinghua.edu.cn.}}
\date{}
\maketitle

\begin{abstract} This paper is an extension of
\cite{K-L}. In this paper,
we will study the blowup behavior
of a surface sequence $\Sigma_k$ immersed in $\R^n$ with
bounded Willmore functional and fixed genus $g$.
We will prove that, we can
decompose $\Sigma_k$ into finitely many  parts:
$$\Sigma_k=\bigcup_{i=1}^m\Sigma_k^i,$$
and find $p_k^i\in \Sigma_k^i$,  $\lambda_k^i
\in\R$, such that $\frac{\Sigma_k^i-p_k^i}
{\lambda_k^i}$
converges locally in the sense of
varifolds to a complete branched immersed surface
$\Sigma_\infty^i$ with
$$\sum_i\int_{\Sigma_\infty^i}K_{\Sigma_\infty^i}=2\pi(2-2g).$$
The basic tool we use in this paper is a generalized convergence theorem
of F. H\'elein.
\end{abstract}



\date{}
\maketitle

\section{Introduction}
For an immersed surface $\ f : \Sigma
\rightarrow \R^n\ $ the Willmore functional
is defined by
\begin{displaymath}
    W(f) = \frac{1}{4} \int_\Sigma |H_f|^2 d \mu_{f},
\end{displaymath}
where $H_f=\Delta_{g_f}f$ denotes the mean curvature vector of $f$, $g_f = f^*
g_{euc}$ the pull-back metric and $\mu_f$ the induced area measure
on $\Sigma$.
This functional first appeared in the papers of Blaschke \cite{Bl}
and Thomsen \cite{T}, and was reinvented and popularized
by Willmore \cite{W}.

We denote the infimum of Willmore functional of immersed surfaces
of genus $p$ by $\beta_p^n$. We have $\beta_p^n\geq 4\pi$ by
Gauss-Bonnet formula, and $\beta_p^n<8\pi$ as observed by Pinkall
and Kusner \cite{K} independently. Willmore conjectured
that $\beta_1^n$ is attained by Clifford torus. This conjecture is
still open.

Given a surface sequence with bounded Willmore functional
and measure, we are  particularly interested to know
what  the limit looks like?
In other words, we expect to  understand
the blowup behavior of such a surface sequence.
It is very important
as we meet blowup almost
everywhere in the
study of Willmore functional.
For example,  if   $\Sigma_t$
is a Willmore flow defined on $[0,T)$, then
by $\epsilon$-regularity proved
in \cite{K-S},  $\int_{B_\rho\cap
\Sigma_t}|A_t|^2<\epsilon$ implies
$\|\nabla_{g_t}^mA_t\|_{L^\infty(B_\frac{\rho}{2}
\cap \Sigma_t)}
<C(m,\rho)$.  Then
 $\Sigma_t$  converges smoothly
in any compact subset of $\R^n$ minus the concentration points set which is defined by
$$\mathcal{S}=\{p\in\R^n:\lim_{r\rightarrow 0}
\liminf_{t\rightarrow T}
\int_{B_r(p)\cap\Sigma_t} |A_t|^2>0\}.$$
So, if we want to have a good knowledge of Willmore flow, we
have to learn  the
behavior, especially the structure of the bubble trees
of $\Sigma_t$ near the concentration
points.

Note that $W(f_k)<C$ implies
$\int_{\Sigma}|A_{k}|^2
d\mu_k<C'$.
One expects that
 $\|f_k\|_{W^{2,2}}$ is equivalent to
$\int|A_k|^2d\mu_k=\int g_k^{ij}g_k^{km}A_{ik}A_{jm}
\sqrt{|g_k|}dx$. However, it is not always true.
One reason is that  the diffeomorphism group
of a surface is extremely big. Therefore, even when an
immersion sequence $f_k$ converges smoothly,
we can easily find a diffeomorphism sequence
$\phi_k$ such that $f_k\circ \phi_k$ will not converge.
Moreover,
the Sobolev embedding
 $W^{2,2q}\hookrightarrow C^1$ is invalid when $q=1$,
so that it is impossible to  estimate the $L^\infty$
norms of
$g^{-1}_k$ and $g_{k}$ via the Sobolev inequalities
 directly.

To overcome these difficulties, an
approximate decomposition lemma was used by L. Simon
when he proved the existence of the minimizer
\cite{S}. He proved that $\beta_p^n$ can be attained if
$p=1$ or
\begin{equation}\label{simon}
p>1,\s and\s \beta_p^n< \omega_p^n=\min\Big\{4\pi+\sum\limits_{i}(\beta_{p_i}^n-4\pi):
\sum\limits_{i} p_i =p,\,1 \leq p_i < p\Big\}.
\end{equation}
Then Bauer and Kuwert proved that \eqref{simon} is
always true,
thus $\beta_p^n$ can be attained for any $p$ and $n$ \cite{B-K}.
Later, such a technique was extended by
W. Minicozzi to get the
minimizer of $W$ on Lagrangian tori \cite{M},
by Kuwert-Sch\"atzle to get
the  minimizer of $W$ in a fixed conformal
class \cite{K-S3}, and by  Sch\"atzle
to get the minimizer of $W$ with boundary
condition \cite{Sh}.

In a recent paper \cite{K-L}, we presented a new approach.
Given an  immersion sequence
$f_k$, we consider each $f_k$ as a conformal immersion
of $(\Sigma,h_k)$ in $\R^n$, where $h_k$ is the
smooth metric with  Gaussian curvature $\pm1$ or 0.
On the one hand, the conformal
diffeomorphism group of $(\Sigma,h_k)$ is very small.
 On
the other hand, if we set $g_{f_k}=e^{2u_k} g_{euc}$ on an
isothermal coordinate system, then we can estimate
$\|u_k\|_{L^\infty}$ from the compensated compactness property  of
$K_{f_k}e^{2u_k}$. Thus we may  get the upper boundary of
$\|f_k\|_{W^{2,2}}$ via the equation $\Delta_{h_k}f_k=H_{f_k}$.
However, the compensated compactness only holds when the $L^2$ norm
of the second fundamental form is small locally, thus the blowup
analysis is needed here. Our basic tools  are the following 2
results:
\begin{thm}\cite{H}\label{Helein} Let $f_k\in W^{2,2}(D,\R^n)$
be a sequence of conformal immersions with induced metrics
$(g_k)_{ij} = e^{2u_k} \delta_{ij}$, and assume
$$
\int_D |A_{f_k}|^2\,d\mu_{g_k} \leq \gamma <
\gamma_n =
\begin{cases}
8\pi & \mbox{ for } n = 3,\\
4\pi & \mbox{ for }n \geq 4.
\end{cases}
$$
Assume also that $\mu_{g_k}(D) \leq C$ and $f_k(0) = 0$.
Then $f_k$ is bounded in $W^{2,2}_{loc}(D,\R^n)$, and there
is a subsequence such that one of the following two alternatives
holds:
\begin{itemize}
\item[{\rm (a)}] $u_k$ is bounded in $L^\infty_{loc}(D)$ and
$f_k$ converges weakly in $W^{2,2}_{loc}(D,\R^n)$ to a conformal
immersion $f \in W^{2,2}_{loc}(D,\R^n)$.
\item[{\rm (b)}] $u_k \to - \infty$ and $f_k \to 0$ locally uniformly on $D$.
\end{itemize}
\end{thm}

\begin{thm}\label{D.K.}\cite{D-K}
Let $h_k,h_0$ be smooth Riemannian metrics on a surface $M$,
such that $h_k \to h_0$ in $C^{s,\alpha}(M)$, where $s \in \N$,
$\alpha \in (0,1)$. Then for each $p \in M$ there exist
neighborhoods $U_k, U_0$ and smooth conformal diffeomorphisms
$\vartheta_k:D \to U_k$, such that $\vartheta_k \to \vartheta_0$
in $C^{s+1,\alpha}(\overline{D},M)$.
\end{thm}
A $W^{2,2}$-conformal
immersion is defined as follows:
\begin{defi}\label{defconformalimmersion}
Let $(\Sigma,g)$ be a Riemann surface. A map $f\in W^{2,2}(\Sigma,g,\mathbb{R}^n)$
is called a conformal immersion, if the induced metric
$g_{f} = df\otimes df$ is given by
$$
g_{f} = e^{2u} g \quad \mbox{ where } u \in L^\infty(\Sigma).
$$
For a Riemann surface $\Sigma$ the set of all
$W^{2,2}$-conformal immersions  is denoted by
$\CI(\Sigma,g,\R^n)$. When $f\in W^{2,2}_{loc}
(\Sigma,g,\R^n)$ and $u\in L^\infty_{loc}(\Sigma)$, we say
$f\in W^{2,2}_{conf,loc}(\Sigma,g,\R^n)$.
\end{defi}

\begin{rem}
F. H\'elein first proved Theorem \ref{Helein} is true for $\gamma<
\frac{8\pi}{3}$ \cite[Theorem 5.1.1]{H}. In \cite{K-L}, we  show
that the constant $\gamma_n$ is optimal.
\end{rem}
\noindent Theorem \ref{Helein} together with Theorem \ref{D.K.}
give the convergence of a $W^{2,2}$-conformal sequence
of $(D,h_k)$ in $\R^n$ with $h_k$ converging smoothly
to $h_0$.

Then using the theory of moduli space of Riemann surface,
we proved in \cite{K-L} the following

\begin{thm}\label{KL}\cite{K-L} Let $f\in W^{2,2}_{conf}
(\Sigma,h_k,\R^n)$. If
\begin{equation}\label{omega}
W(f_k)\leq \left\{\begin{array}{ll}
                    8\pi-\delta&p=1\\
                    \min\{8\pi,\omega_p\}-\delta&p>1
                 \end{array}\right.,\s \delta>0,
\end{equation}
then
the conformal class sequence
represented by $h_k$ converges in $\mathcal{M}_p$.
\end{thm}
In other words, $h_k$
converges to a metric $h_0$ smoothly. This was
also proved by T. Rivi\`{e}re \cite{R}. Then up to
M\"obius transformations, $f_k$
will converge weakly in $W^{2,2}_{loc}(\Sigma
\setminus\{\mbox{finite points}\},h_0)$ to a
$W^{2,2}(\Sigma,h_0)$-conformal immersion. In this way,
we give a new proof of the existence of
minimizer of Willmore functional with fixed genus.

\eqref{omega} also gives us a hint that,
it is the degeneration of complex structure that makes
the trouble for the convergence of an immersion sequence
with
\begin{equation}\label{bmw}
\mu(f_k)+W(f_k)<C.
\end{equation}

In \cite{C-L}, the Hausdorff limit of $\{f_k\}$
with \eqref{bmw} was studied,
using conformal immersion as a tool.
We proved that, the limit of $f_0$ is a conformal
branched immersion from a stratified
surface $\Sigma_\infty$  into $\R^n$.
Briefly speaking,
 if $(\Sigma_0,h_0)$ is the limit of $(\Sigma,h_k)$
in $\overline{\mathcal{M}_p}$, then
$f_k$ converges weakly in the $W^{2,2}$ sense
in any component of $\Sigma_0$ away from the blowup points
$$\mathcal{S}(f_k)=\{p\in D:\lim_{r\rightarrow 0}
\liminf_{k\rightarrow+\infty}\int_{B_r(p,h_0)}|A(f_k)|^2d\mu_{f_k}\geq 4\pi\}.$$
Meanwhile, some  bubble trees, which consist of $W^{2,2}$ branched conformal
immersions of $S^2$ in $\R^n$ will appear.
As a corollary, we get the following
\begin{pro}\cite{C-L}
Let  $f_k:\Sigma\rightarrow \R^n$ be a sequence of smooth
immersions with \eqref{bmw}.
Assume the Hausdorff limit of $f_k(\Sigma)$
is not a union of $W^{2,2}$ branched conformal immersed
spheres. Then the  complex structure of $c_k$ induced by $f_k$ diverges in
the moduli space if and only if  there are a seqence of closed
curves $\gamma_k$
which are nontrivial in $H^1(\Sigma)$, such that
the length of $f_k(\gamma_k)$ converges to 0.
\end{pro}

Thus, when the conformal class induced by $f_k$
diverges in the moduli space,  topology will
be lost. They are two reasons why the topology is lost.
One reason is that Theorem \ref{Helein} does not ensure
the limit is an immersion on each component of $\Sigma_0$.
If $f_k$ converges to a point in some components, then
some topologies are taken away. The other reason is that on each collar which is conformal
to $Q(T_k)=S^1\times[-T_k,T_k]$ with $T_k\rightarrow+\infty$,
there must exist a sequence
$t_k\in[-T_k,T_k]$ such that $f_k(S^1\times\{t_k\})$  will shrink to a point.

It is not easy to calculate  how many topologies are lost, but it
is indeed possible to find where $\int_\Sigma K_{f_k}d\mu_{f_k}$
is lost. We have to study those bubbles
which have nontrivial topologies but shrink to points.
For this sake, we should check if those conformal
immersion sequences which converge to points
will converge to   immersions after being rescaled:

\begin{thm}\label{convergence}
Let $\Sigma$ be a smooth connected Riemann surface without boundary,
and $\Omega_k\subset\subset\Sigma$ be domains with
$$\Omega_1\subset \Omega_2\subset\cdots\Omega_k\subset\cdots,\s
\bigcup_{i=1}^\infty\Omega_i=\Sigma.$$
Let $\{h_k\}$ be a  smooth metric sequence over $\Sigma$
which converges to $h_0$ in $C^\infty_{loc}(\Sigma)$, and
$\{f_k\}$ be
a conformal immersion sequence of $(\Omega_k,h_k)$
in $\R^n$ satisfying
\begin{itemize}
\item[{\rm 1)}] $\mathcal{S}(f_k):=
\{p\in\Sigma: \lim\limits_{r\rightarrow 0}\liminf\limits_{k\rightarrow+\infty}
\int_{B_r(p,h_0)}|A_{f_k}|^2d\mu_{f_k}\geq 4\pi \}=\emptyset$.
\item[{\rm 2)}] $f_k(\Omega_k)$ can be extended to a closed compact
immersed surface $\Sigma_k$ with
$$\int_{\Sigma_k}(1+|A_{f_k}|^2)d\mu_{f_k}<\Lambda.$$
\end{itemize}
Take a curve $\gamma:[0,1]\rightarrow \Sigma$,
and set $\lambda_k=diam\, f_k(\gamma[0,1])$.
Then
we can find  a subsequence of $\frac{f_k
-f_k(\gamma(0))}{\lambda_k}$ which
converges  weakly in
$W^{2,2}_{loc}(\Sigma)$ to an
$f_0\in W_{conf,loc}^{2,2}(\Sigma,\R^n)$.
Further, we can find an inverse $I=\frac{y-y_0}{|y-y_0|^2}$
with $y_0\notin f_0(\Sigma)$ such that
$$\int_\Sigma(1+|A_{I(f_0)}|^2)d\mu_{I(f_0)}<+\infty.$$
\end{thm}

When $\Sigma$ is a compact closed surface minus finitely
many points,
$f_0$ may not be compact. However,
by Removability of singularity (see Theorem \ref{removal} in
section 2), $I(f_0)$ is a conformal branched immersion.
Thus $f_0$ is complete.

\begin{defi}
We  call $f$   a generalized  limit of $f_k$, if we can
find a point $x_0\notin \mathcal{S}(f_k)$ and a positive
sequence $\lambda_k$ which is equivalent to 1 or tends to 0,
such that
$\frac{f_k-f_k(x_0)}{\lambda_k}$ converges to $f$ weakly in
$W^{2,2}_{loc}(\Sigma\setminus \mathcal{S}(f_k))$.
\end{defi}

Obviously, if $f$ and $f'$ are both generalized limits
of $f_k$, then $f=\lambda f'+b$ for some $\lambda$ and $b$.
We will not distinguish between $f$ and $f'$.

Near the concentration points, we will get
some bubbles. The divergence of complex structure also
gives us some bubbles.  In
\cite{C-L}, we  only considered the bubbles with
$\lambda_k\equiv1$.
In this paper, we will study  the bubbles with $\lambda_k
\rightarrow 0$ which do not appear in the
Hausdorff limit. All the bubbles
can be considered as  conformal branched
immersions from $\C$ (or $S^1\times \R$, $S^2$) into
$\R^n$.
 However,
the structures of bubble trees here
are  much more complicated than those of
harmonic maps. For example, there might
exist infinite many bubbles here, therefore, we should
neglect  the bubbles  which do not carry
total Gauss curvature.

\begin{defi}
We say a conformal branched immersion of $S^1\times\R$
into $\R^n$ is trivial, if for any $t$,
$$\int_{S^1\times\{t\}}\kappa\neq 2m\pi+\pi,\s
for\s some\s m\in\mathbb{Z}.$$
\end{defi}

The bubble trees constructed in this paper
consist of finitely many branches.
Small branches are on the big branches level by level.
Each branch consists of nontrivial bubbles, bubbles
with concentration, and the first bubble (see definitions in
Section 4).
We can classify the bubbles into four types:
$T_\infty$, $T_0$, $B_\infty$ and $B_0$
(see Definition \ref{typeofbubble}). We will show that a $T_0$ type bubble must follow a $B_\infty$ type bubble,
and a $T_\infty$ type bubble must follow a $B_0$ type
bubble.

Moreover, we have total Gauss curvature identity.
To state the total Gauss curvature identity precisely, we
have to divide it into 3 cases.\vspace{0.7ex}

{\bf Hyperbolic case (genus$>1$):} Let $\Sigma_0$ be the stable surface in $
\overline{\mathcal{M}}_g$ with nodal points $\mathcal{N}=\{
a_1,\cdots, a_{m}\}$.
$\Sigma_0$ is obtained by pinching
some curves in a surface to points,
thus $\Sigma_0\setminus\mathcal{N}$ can be divided
into finitely many components $\Sigma_0^1$, $\cdots$,
$\Sigma_0^s$. For each $\Sigma_0^i$, we can
extend $\Sigma_0^i$ to a smooth closed Riemann surface $\overline{\Sigma_0^i}$
by adding a point at each puncture. Moreover, the
complex structure of $\Sigma_0^i$ can be extended
smoothly to a complex structure of $\overline{\Sigma_0^i}$.

We say $h_0$ to be a hyperbolic structure on $\Sigma_0$ if $h_0$ is a  smooth complete metric on
$\Sigma_0\setminus\mathcal{N}$ with finite volume and Gauss curvature $-1$.
We define $\Sigma_{0}(a_j,\delta)$ to be the domain in $\Sigma_0$ which satisfies
$$a_j\in \Sigma_0(a_j,\delta),\s and \s injrad_{\Sigma_0\setminus\mathcal{N}}^{h_0}(p)<\delta\s \forall p\in\Sigma_0(a_j,\delta)\setminus\{a_j\}.$$
We set $h_0^i$ to be a
smooth metric over $\overline{\Sigma_0^i}$
which
is conformal to $h_0$ on $\Sigma_0^i$. We may assume $h_0^i$
has curvature $\pm1$ or curvature $0$ and measure 1.

Now, we let $\Sigma_k$ be a sequence of compact Riemann
surfaces of
fixed genus $g$ whose metrics $h_k$ have curvature $-1$,
such that $\Sigma_k
\rightarrow \Sigma_0$ in
$\overline{\mathcal{M}_g}$.
Then, there exist
a maximal collection $\Gamma_k = \{\gamma_k^1,\ldots,\gamma_k^{m}\}$
of pairwise disjoint, simply closed geodesics in $\Sigma_k$
with $\ell^j_k = L(\gamma_k^j) \to 0$, such that after passing
to a subsequence the following hold:
\begin{itemize}
\item[{\rm (1)}] There are maps $\varphi_k \in C^0(\Sigma_k,\Sigma_0)$,
such that $\varphi_k: \Sigma_k \backslash \Gamma_k \to \Sigma_0 \backslash \mathcal{N}$
is diffeomorphic and $\varphi_k(\gamma_k^j) = a_j$ for $j = 1,\ldots,m$.
\item[{\rm (2)}] For the inverse diffeomorphisms
$\psi_k:\Sigma_0 \backslash \mathcal{N} \to \Sigma_k \backslash \Gamma_k$,
we have $\psi_k^\ast (h_k) \to h_0$ in $C^\infty_{loc}(\Sigma_0 \backslash
\mathcal{N})$.
\item[{\rm (3)}] Let $c_k$ be the complex structure over $\Sigma_k$, and $c_0$
be the complex structure over $\Sigma_0\setminus\mathcal{N}$. Then
$$\psi_{k}^*(c_k)\rightarrow c_0\s in\s
C^\infty_{loc}(\Sigma_0\setminus\mathcal{N}).$$
\item[{\rm (4)}]For each $\gamma_k^j$, there is a
collar $U_k^j$ containing $\gamma_k^j$, which is
isometric to cylinder
$$Q_k^j=S^1\times(-\frac{\pi^2}{l_k^j},\frac{\pi^2}{l_k^j}),\s
with\s metric\s
h_k^j=\left(\frac{1}{2\pi\sin(\frac{l_k^j}{2\pi}t+\theta_k)}\right)^2(dt^2+d\theta^2),$$
where $\theta_k=\arctan(\sinh(
\frac{l_k^j}{2}))+\frac{\pi}{2}$.
Moreover, for any $(\theta,t)\in S^1\times
(-\frac{\pi^2}{l_k^j},\frac{\pi^2}{l_k^j})$, we have
\begin{equation}\label{injrad}
\sinh(injrad_{\Sigma_k}(t,\theta))\sin(
\frac{l_k^jt}{2\pi}+\theta_k)
=\sinh\frac{l_k^j}{2}.
\end{equation}
Let $\phi_k^j$ be the isometric between $Q_k^j$ and $U_k^j$. Then
$\varphi_k\circ\phi_k^{j}(T_k^j+t,\theta)\cup
\varphi_k\circ\phi_k^{j}(-T_k^j+t,\theta)$ converges in $C^\infty_{loc}((-\infty,0)\cup (0,\infty))$ to an isometric from $S^1\times(-\infty,0)\cup S^1\times(0,+\infty)$ to $\Sigma_0(a_j,1)\setminus \{a_j\}$.

\end{itemize}

Items 1) and 2) in the above  can be found in Proposition 5.1 in \cite{Hum}.
The main part of 3) is just the collar Lemma.

Now, we consider a sequence $f_k\in W^{2,2}_{conf}
(\Sigma,h_k,\R^n)$, with
$$\mu(f_k)+W(f_k)<\Lambda.$$
By Theorem \ref{convergence}, on each component
$\Sigma_k^i$, $f_k\circ \psi_k$ has a generalized limit
$f_0^i\in W^{2,2}_{conf}(\overline{\Sigma_k^i}\setminus
A^i,h_0^i,\R^n)$,
where $A^i$ is a finite set. We have the following

\begin{thm}\label{main}
Let $f^1$, $f^2$, $\cdots$ be
all of the non-trivial bubbles of $\{f_k\}$. Then
$$\sum_i\int_{\overline{\Sigma_k^i}}K_{f_0^i}d\mu_{f_0^i}+
\sum_i\int_{S^2}K_{f^i}d\mu_{\varphi^i}=2\pi\chi(\Sigma).$$
\end{thm}\vspace{0.7ex}

{\bf Torus case:} Let $(\Sigma,h_k)=\C/(\pi,z)$, where
$|z|\geq\pi$ and $|\Re{z}|\leq\frac{\pi}{2}$.
We can write
$$(\Sigma,h_k)=S^1\times\R/G_k,$$
where $S^1$ is the circle with perimeter 1 and
$G_k\cong \mathbb{Z}$ is the transformation group generalized by
$$(t,\theta)\rightarrow (t+a_k,\theta+\theta_k),\s where\s
a_k\geq \sqrt{\pi^2-\theta_k^2},\s and\s \theta_k\in [-\frac{\pi}{2},\frac{\pi}{2}].$$
$(\Sigma_k,h_k)$ diverges
in $\mathcal{M}_1$ if and only if $a_k\rightarrow+\infty$.

Then any $f_k\in W^{2,2}_{conf}(\Sigma,h_k,\R^n)$ can be
lifted  to a conformal immersion  $f_k':S^1\times\R
\rightarrow\R^n$ with
$$f_k'(t,\theta)=f_k'(t+a_k,\theta+\theta_k).$$
After translating, we may assume that
$f_k'(-t+\frac{a_k}{2},\theta)$ and $f_k'(t-\frac{a_k}{2},\theta)$
have no concentrations. We let $\lambda_k=diam
f_k'(S^1\times{\frac{a_k}{2}})$,
then $\frac{f_k'(-t+\frac{a_k}{2},\theta)-f(\frac{a_k}
{2},0)}{\lambda_k}$ and $\frac{f_k'(t-\frac{a_k}{2},\theta)
-f_k'(\frac{a_k}{2},\theta_k)}{\lambda_k}$ will
converge to $f_0^1$ and $f_0^2$ respectively in
$W^{2,2}_{loc}(S^1\times[0,+\infty))$. However, they can be
glued together via
$$f_0=\left\{\begin{array}{ll}
                f_0^1(-t,\theta)&t\leq 0\\
                f_0^2(t,\theta+\theta_0)&t>0,
            \end{array}\right.$$
into a conformal immersion of $S^1\times\R$
in $\R^n$,
where $\theta_0=\lim\limits_{k\rightarrow+\infty}\theta_k$.
Then we have
\begin{thm}\label{main2}
$$\int_{S^1\times\R}K_{f_0}d\mu_{f_0}
+\sum_{i=1}^m\int_{S^1\times\R}K_{f^i}d\mu_i=0,$$
where $f^1$, $\cdots$, $f^m$ are all of the non-trivial bubbles of $f_k'$.
\end{thm}\vspace{0.7ex}

{\bf Sphere case:} When $\Sigma$ is the sphere, we can
let $h_k\equiv h_0$. There is no bubble from collars.
We have
\begin{thm}\label{sphere} Let $f_0$ be the generalized limit of $f_k$.
Then
$$\int_{S^2}K_{f_0}d\mu_{f_0}+\sum_{i=1}^m\int_{S^1\times\R}
K_{f^i}d\mu_{f^i}=4\pi,$$
where $f^1$, $\cdots$, $f^m$ are all of the  non-trivial bubbles.
\end{thm}

Put Theorem \ref{main}--\ref{sphere} together, we get
the main theorem of this paper, which is
a precise version of Theorem \ref{KL}:
\begin{thm}
Let $\Sigma_k$ be a sequence of surfaces
immersed in
$\R^n$ with bounded Willmore functional.
Assume $g(\Sigma_k)=g$. Then we can
decompose $\Sigma_k$ into finite parts:
$$\Sigma_k=\bigcup_{i=1}^m\Sigma_k^i,\s \Sigma_i\cap\Sigma_j
=\emptyset,$$
and find $p_k^i\in \Sigma_k^i$,  $\lambda_k^i
\in\R$, such that $\frac{\Sigma_k^i-p_k^i}
{\lambda_k^i}$
converges locally in the sense of
varifolds to a complete branched immersed surface
$\Sigma_\infty^i$ with
$$\sum_i\int_{\Sigma_\infty^i}K_{\Sigma_\infty^i}=2\pi(2-2g),
\s and\s \sum_{i}W(\Sigma_\infty^i)\leq \lim_{k\rightarrow+\infty}
W(\Sigma_k).$$

\end{thm}


\begin{rem}
Parts of Theorem \ref{sphere} have appeared in \cite{L-L},
in which we assumed that
$\{f_k\}\subset W^{2,2}_{conf}(D,\R^n)$
and does not converge to a point.
\end{rem}

\section{Preliminary}

\subsection{Hardy estimate}
Let $f\in W^{2,2}_{conf}(D,\R^n)$ with $g_f=
e^{2u}(dx^1\otimes dx^1+dx^2\otimes dx^2)$ and
$\int_D|A_f|^2<4\pi-\delta$. $f$ induces a
Gauss map
$$G(f)=e^{-2u}(f_1\wedge f_2):D\rightarrow G(2,n)\hookrightarrow
\C\mathbb{P}^{n-1}.$$
Following \cite{M-S}, we define
 the map $\Phi(f):\C \to \C P^{n-1}$ by
$$
\Phi(f)(z) =
\left\{\begin{array}{ll}
G(f)(z) & \mbox{ if }z \in D\\
G(f)(\frac{1}{\overline{z}}) & \mbox{ if }z \in \C \backslash \overline{D}.
\end{array}\right.$$
Then $\Phi(f)\in W_0^{1,2}(\mathbb{C},\C P^{n-1})$ and
$\int_{\mathbb{C}}{\Phi}^*(f) (\omega) =0$, where $\omega$
is the K\"ahler form of $\C\mathbb{P}^{n-1}$. Thus by  Corollary 3.5.7
in \cite{M-S},
$\Psi(f)=*\Phi^*(f)(\omega)$ is in Hardy space, and
\begin{equation}\label{Psi}
\|\Psi(f)\|_{
\mathcal{H}}<C(\delta)\|A_f\|_{L^2(D)}.
\end{equation}
Note that
\begin{equation}\label{psi-k}
\Psi(f)|_{D}=K_{f}e^{2u}.
\end{equation}

If we set that $v$  solves the equation $-\Delta v=\Psi(f)$,
$v(\infty)=0$, then we have
\begin{equation*}
\|v\|_{L^\infty(\R^n)}+\|\nabla v\|_{L^2(\R^n)}+\|\nabla^2 v\|_{L^1(\R^n)}<
C\|\Psi(f)\|_{\mathcal{H}}.
\end{equation*}
Noting that $u-v$ is harmonic on $D$, we get
\begin{equation}\label{hardy0}
\|u\|_{L^\infty(D_\frac{1}{2})}+\|\nabla u\|_{L^2(D_\frac{1}{2})}+\|\nabla^2 u\|_{L^1(D_\frac{1}{2})}<
C(\|\Psi(f)\|_{\mathcal{H}}+\|u\|_{L^1(D)}).
\end{equation}

\subsection{Gauss-Bonnet formula}
Let $f\in W^{2,2}_{conf}(\Sigma,g,\R^n)$ with
$g_f=e^{2u}g$.  Let $\gamma$ be
a smooth curve.
On $\gamma$, we define
\begin{equation}\label{geodesic.curvature}
\kappa_{f}=\frac{\partial u}{\partial n}+\kappa_g,
\end{equation}
where $n$ is one of the
unit normal field along $\gamma$ which is compatible to
$\kappa_g$. By \eqref{hardy0}, $\frac{\partial u}{\partial n}$
is well-defined.
In \cite{K-L}, we proved that $u$ satisfies the weak equation
$$-\Delta_g u=K_{f}e^{2u}-K_g.$$
Then, for any domain $\Omega$ with smooth boundary,
we have the Gauss-Bonnet formula:
$$\int_{\partial\Omega}\kappa_f=\chi(\overline{\Omega})+\int_{\Omega}
K_fd\mu_f.$$

\subsection{Convergence of $\int K_{f_k}d\mu_{f_k}$}

By \eqref{Psi} \eqref{psi-k} \eqref{hardy0} and
Theorem \ref{Helein}, we have :

\begin{lem}\label{measureconvergence}
Let $f_k$ be a conformal sequence from $D$ into $\R^n$
with $g_{f_k}=e^{2u_k}g_0$ and
$\int_D|A_{f_k}|^2d\mu_f\leq \gamma<4\pi$, which converges
to $f_0$ weakly. We assume $f_0$ is not a point map, and
$g_{f_0}=e^{2u_0}g_0$.
Then we can find a subsequence, such that
\begin{equation}\label{hardy}
K_{f_k}d\mu_{f_k} \rightharpoonup
K_{f_0}d\mu_{f_0}\s over \s D_\frac{1}{2},\s \mbox{in distribution,}
\end{equation}
and
\begin{equation*}
u_k\rightharpoonup u_0,\s in\s W^{1,2}(D_\frac{1}{2}).
\end{equation*}
\end{lem}

We will use the following
\begin{cor}
Let  $f_k$ be a conformal sequence of $D\setminus
D_\frac{1}{2}$ in $\R^n$, which converges to $f_0\in W^{2,2}_{{conf},loc}(D\setminus D_\frac{1}{2},\R^n)$.
For any $t\in(\frac{1}{2},1)$ with $\partial D_t\cap \mathcal{S}(f_k)
=\emptyset$, we have
$$\lim_{k\rightarrow+\infty}\int_{\partial D_t}\kappa_{f_k} ds_k=
\int_{\partial D_t}\kappa_{f_0}ds_0.$$
\end{cor}

\proof Take $s\in (t,1)$, such that $\mathcal{S}(f_k)\cap
\overline{D_s\setminus D_t}=\emptyset$.
Let $g_{f_k}=e^{2u_k}g_0$ and $\varphi\in C^\infty_0(D_s)$,
which is 1 on $D_t$.
Then we have
$$-\int_{\partial D_t}\frac{\partial u_k}{\partial r}ds=-
\int_{D_s\setminus D_t}\nabla u_k\nabla\varphi
d\sigma
+\int_{D_s\setminus D_t}\varphi K_ke^{2u_k}d\mu_{f_k},$$
and the right-hand side will converge to
$$-
\int_{D_s\setminus D_t}\nabla u_0\nabla\varphi
d\sigma+\int_{D_s\setminus D_t}\varphi K_0e^{2u_0}d\mu_{f_0},\s
as\s k\rightarrow+\infty.$$
Then we get
$$-\int_{\partial D_t}\frac{\partial u_k}{\partial r}ds
\rightarrow -\int_{\partial D_t}\frac{\partial u_0}{\partial r}ds.$$
By \eqref{geodesic.curvature} we
get
$$\int_{\partial D_t}\kappa_k\rightarrow
\int_{\partial D_t}\kappa_0.$$
\endproof

\subsection{Removability of singularity}
We have the following
\begin{thm}\label{removal}\cite{K-L}
Suppose that $f\in W^{2,2}_{{conf},loc}(D\backslash \{0\},\R^n)$ satisfies
$$
\int_D |A_f|^2\,d\mu_g < \infty \quad \mbox{ and } \quad \mu_g(D) < \infty,
$$
where $g_{ij} = e^{2u} \delta_{ij}$ is the induced metric. Then
$f \in W^{2,2}(D,\R^n)$ and we have
\begin{eqnarray*}
u(z) & = & m\log |z|+ \omega(z) \quad \mbox{ where }
m\geq 0,\, z\in \mathbb{Z},\,\omega \in C^0 \cap W^{1,2}(D),\\
-\Delta u & = & -2m\pi \delta_0+K_g e^{2u} \quad \mbox{ in }D.
\end{eqnarray*}
The multiplicity of the immersion at $f(0)$ is given by
$$
\theta^2\big(f(\mu_g \llcorner D_\sigma(0)),f(0)\big) = m+1 \quad \mbox{ for any small }
\sigma > 0.
$$
Moreover, we have
\begin{equation}\label{kappa2}
\lim_{t\rightarrow 0}\int_{\partial D_t}\kappa_{f}
ds_f=2\pi (m+1).
\end{equation}
\end{thm}

\proof  We only prove \eqref{kappa2}. For the
proof of other part of the theorem, one can refer to
\cite{K-L}.

Observe that
$$|\int_{\partial D_t}\frac{\partial u}
{\partial r}-\int_{\partial D_{t'}}\frac{\partial u}
{\partial r}|=|\int_{D_t\setminus D_{t'}}
Kd\mu|\rightarrow 0$$
as $t, t'\rightarrow 0$. Then
$\lim\limits_{t\rightarrow 0}\int_{\partial D_t}\frac{\partial u}
{\partial r}$ exists.

Since $\omega\in W^{1,2}(D_r)$,
we can find $t_k\in [2^{-k-1},2^{-k}]$, s.t.
$$(2^{-k}-2^{-k-1})\int_{\partial D_{t_k}}|
\frac{\partial w}{\partial r}|
=\int_{2^{-k-1}}^{2^{-k}}(\int_{\partial D_t}
|\frac{\partial w}{\partial r}|)dt\leq C\|\nabla w\|_{L^2
(D_{2^{-k}})}2^{-k},$$
which implies that
$\int_{\partial D_{t_k}}\frac{\partial w}{\partial r}
\rightarrow 0$. Then we get
$\int_{\partial D_{t_k}}\frac{\partial u}
{\partial r}\rightarrow 2\pi m$, which implies
that
$$\lim_{t\rightarrow 0}
\int_{\partial D_{t}}\frac{\partial u}
{\partial r}\rightarrow 2\pi m.$$

\endproof

\begin{rem} In the proof of Theorem \ref{removal} in \cite{K-L}, we get that
\begin{equation}\label{isolate}
\lim_{z\rightarrow 0}\frac{|f(z)-f(0)|}{|z|^{m+1}}=\frac{e^{w(0)}}{m+1}.
\end{equation}
\end{rem}

We give the following definition:
\begin{defi} \label{defconformalimmersion}
A map $f\in W^{2,2}(\Sigma,\mathbb{R}^n)$
is called a $W^{2,2}$- branched  conformal immersion, if we can find
finitely many points $p_1$, $\cdots$, $p_m$, s.t.
$f\in W^{2,2}_{conf,loc}(\Sigma\setminus\{p_1,\cdots,p_m\})$, and
$$
\mu(f)<+\infty,\s \int_{\Sigma}|A_f|^2d\mu_f<+\infty.
$$
\end{defi}

For the behavior at infinity of complete conformally
parameterized surfaces,  we have the following

\begin{thm}\label{removal2}
Suppose that $f\in W^{2,2}_{{conf},loc}(\C\setminus D_R,\R^n)$ with
$$
\int_{\C\setminus D_R} |A_f|^2\,d\mu_{g} < \infty,
$$
where $g_{ij} = e^{2u} \delta_{ij}$ is the induced
metric. We assume $f(\C\setminus D_{2R})$
is complete. Then we have
\begin{equation*}
u(z) =  m\log |z|+ \omega(z) \quad
\mbox{ where }
m\geq 0,\, z\in \mathbb{Z},\,\omega \in  W^{1,2}(\C\setminus D_{2R}).
\end{equation*}
Moreover, we have
\begin{equation}\label{kappa3}
\lim_{t\rightarrow +\infty}\int_{\partial D_t}\kappa_{f}
ds_f=2\pi (m+1).
\end{equation}

\end{thm}

The proof of \eqref{kappa3} is similar to that of
\eqref{kappa2}. Other part of the proof can
be found in \cite{M-S}.
Though Muller-Sverak's result was stated for smooth
surface,  it is easy to check that
their proof also holds for a $W^{2,2}$ conformal immersion.

\section{Proof of Theorem \ref{convergence}}
We first prove the following
\begin{lem}\label{convergence2} Suppose $(\Sigma,h_k)$ to be smooth Riemann surfaces,
where $h_k$ converges to $h_0$ in $C^\infty_{loc}(\Sigma)$.
Let $\{f_k\}\subset W^{2,2}_{{conf},loc}(\Sigma,h_k,\R^n)$ with
$$\mathcal{S}(f_k)=\{p\in\Sigma: \lim_{r\rightarrow 0}
\liminf_{k\rightarrow+\infty}\int_{B_r(p,h_0)}|A_{f_k}|^2d\mu_{f_k}\geq 4\pi\}
=\emptyset.$$
Then $f_k$ converges in $W^{2,2}_{loc}(\Sigma,h_0,\R^n)$
to a point or an $f_0\in W^{2,2}_{conf,loc}(\Sigma,h_0,\R^n)$.
\end{lem}

\proof Let $g_{f_k}=e^{2u_k}h_k$.
We only need to prove the following statement:
for any $p\in\Sigma$, we
can find a neighborhood $V$ which is independent of
$\{f_k\}$, such that $f_k$ converges weakly to $f_0$
in
$W^{2,2}(V,h_0)$. Moreover,
 $\|u_k\|_{L^\infty(V)}<C$ if and only if $f_0\in W^{2,2}_{conf}
(V,\R^n)$; $u_k\rightarrow-\infty$ uniformly,
if and only if $f_0$ is a point map.

Now we prove this statement: Given a point $p$, we choose
 $U_k$, $U_0$, $\vartheta_k$, $\vartheta_0$
 as in the Theorem \ref{D.K.}.
Set $\vartheta_k^*(h_k)=e^{2v_k}g_0$, where
$g_0=(dx^1)^2+(dx)^2$.
We may assume $v_k\rightarrow v_0$ in $C^\infty_{loc}(D)$.

Let $\hat{f}_k=f_k(\vartheta_k)$ which
is a map from $D$ into $\R^n$. It is easy to check that
$\hat{f}_k\in W^{2,2}_{conf}(D,\R^n)$ and
$g_{f_k}=e^{2u_k+2v_k}g_0$. By Theorem \ref{Helein}, we can assume
that $\hat{f}_k$ converges to $\hat{f}_0$
weakly in $W^{2,2}(D_\frac{3}{4})$. Moreover,
$\hat{f}_0$ is  a point when $u_k+v_k\
\rightarrow-\infty$ uniformly on $D_\frac{3}{4}$,
and a conformal immersion when $\sup_{k}
\|u_k+v_k\|_{L^\infty(D_\frac{3}{4})}<+\infty$.

Let $V=\vartheta_0(D_\frac{1}{2})$.
Since $\vartheta_k$ converges to $\vartheta_0$,
$\vartheta_k^{-1}(V)\subset D_\frac{3}{4}$
for any sufficiently large $k$ and
$f_k=\hat{f}_k(\vartheta_k^{-1})$
converges to $f_0=\hat{f}_0
(\vartheta_0^{-1})$ weakly in $W^{2,2}(V,h_0)$.
Moreover,
 $f_0$ is a conformal immersion when $\|u_k\|_{L^\infty(V)}<C$,
and a point when $u_k\rightarrow-\infty$ uniformly in
$V$.

\endproof

{\it The proof of Theorem \ref{convergence}:} When $f_k$ converges to a conformal immersion
weakly, the result is obvious. Now we assume that
$f_k$ converges to a point. For this case,
 $\lambda_k
\rightarrow 0$.

Put $f_k'=\frac{f_k-f_k(\gamma(0))}{\lambda_k}$,
$\Sigma_k'=\frac{\Sigma_k-f_k(\gamma(0))}{\lambda_k}$.
We have two cases:\vspace{0.7ex}

\noindent Case 1: $diam(f_k')<C$. Letting $\rho$ in
inequality (1.3) in \cite{S} tend to infinity,
we get $\frac{\Sigma_k'\cap B_\sigma(\gamma(0))}{\sigma^2}\leq C$ for any $\sigma>0$,
hence we get $\mu(f_k')<C$ by taking
$\sigma=diam(f_k')$. Then Lemma \ref{convergence2}
shows that $f_k'$ converges weakly
in $W^{2,2}_{loc}(\Sigma,h_0)$. Since
$diam\, f_k'
(\gamma)=1$, the weak limit is not a point.\vspace{0.7ex}

\noindent Case 2: $diam(f_k')\rightarrow +\infty$. We take a point
$y_0\in\R^n$ and a constant $\delta>0$, s.t.
$$B_\delta(y_0)\cap \Sigma_k'=\emptyset,\s \forall k.$$
Let $I=\frac{y-y_0}{|y-y_0|^2}$, and
$$f_k''=I(f_k'),\s \Sigma_k''=I(\Sigma_k').$$
By conformal invariance of Willmore functional
\cite{C,W}, we have
$$\int_{\Sigma''}|A_{\Sigma''}|^2d\mu_{\Sigma''}
=\int_{\Sigma}|A_\Sigma|^2d\mu_{\Sigma}<\Lambda.$$
Since $\Sigma_k''\subset B_\frac{1}{\delta}(0)$, also by (1.3) in \cite{S},
we get $\mu(f_k'')<C$. Thus
$f_k''$ converges weakly in $W^{2,2}_{loc}(\Sigma\setminus
\mathcal{S}(f_k''),h_0)$.

Next, we prove that $f_k''$ will not converge to a point by assumption.
If $f_k''$ converges to a point in
$W^{2,2}_{loc}(\Sigma\setminus \mathcal{S}(f_k''))$,
then the limit must be 0,  for $diam\,(f_k')$
converges to $+\infty$.
By the
definition of $f_k''$, we can find a $\delta_0>0$,
such that $f_k''(\gamma)\cap
B_{\delta_0}(0,h_0)=\emptyset$. Thus for any $p\in \gamma([0,1])
\setminus \mathcal{S}(f_k'')$, $f_k''$ will not converge to $0$. A contradiction.

Then we only need to prove that $f_k'$ converges weakly in
$W^{2,2}_{loc}(\Sigma,h_0,\R^n)$.
Let $f_0''$ be the limit of $f_k''$. By Theorem \ref{removal},
$f_0''$ is a branched immersion of $\Sigma$ in $\R^n$.
Let $\mathcal{S}^*=f_0^{''-1}(\{0\})$.
By \eqref{isolate}, $\mathcal{S}^*$ is isolate.

First, we prove that for any $\Omega\subset\subset\Sigma\setminus
(\mathcal{S}^*\cup\mathcal{S}(\{f_k''\})$, $f_k'$
converges weakly in $W^{2,2}(\Omega,h_0,\R^n)$:
Since $f_0''$ is continuous
on $\bar{\Omega}$, we may assume
$dist(0,f_0''(\Omega))>\delta>0$. Then $dist(0,f_k''(\Omega))>\frac{\delta}{2}$
when $k$ is sufficiently large. Noting that $f_k'
=\frac{f_k''}{|f_k''|^2}+y_0$, we get that $f_k'$ converges weakly in
$W^{2,2}(\Omega,h_0,\R^n)$.

Next, we prove that for each
$p\in \mathcal{S}^*\cup\mathcal{S}(\{f_k''\})$, $f_k'$ also converges in
a neighborhood of $p$.
We use the denotation $U_k$, $U_0$, $\vartheta_k$ and
$\vartheta_0$ with $\theta_k(0)=p$ again.
We only need to prove that $\hat{f}_k'=f_k'(\vartheta_k)$ converges
weakly in $W^{2,2}(D_\frac{1}{2})$.

Let $g_{\hat{f}_k'}=e^{2\hat{u}_k'}(dx^2+dy^2)$.
Since $\hat{f}_k'\in W^{2,2}_{conf}
(D_{4r})$ with $\int_{D_{4r}}|A_{\hat{f}_k'}|^2d\mu_{\hat{f}_k'}<4\pi$ when $r$ is
sufficiently small and $k$ sufficiently large,
by the arguments in subsection 2.1,
we can find a $v_k$ solving the equation
$$-\Delta v_k=K_{\hat{f}_k'}e^{2\hat{u}_k'},\s z\in D_r\s and\s \|v_k\|_{L^\infty(D_r)}<C.$$
Since $f_k'$ converges to a conformal
immersion in $D_{4r}\setminus D_{\frac{1}{4}r}$, by Theorem \ref{Helein},
we may assume that
$\|\hat{u}_k'\|_{L^\infty(D_{2r}\setminus
D_r)}<C$.
 Then
$\hat{u}_k'-v_k$ is a harmonic function with
$\|\hat{u}_k'-v_k\|_{L^\infty(\partial D_{2r}(z))}<C$,
then we get $\|\hat{u}_k'(z)-v_k(z)\|_{L^\infty(D_{2r}(z))}<C$
by the Maximum Principle. Thus, $\|\hat{u}_k'\|_{L^\infty(D_{2r})}<C$,
which implies $\|\nabla f_k'\|_{L^\infty(D_{2r})}<C$.
By the equation $\Delta \hat{f}_k'=e^{2\hat{u}_k'}H_{\hat{f}_k'}$, and
the fact that $\|e^{2\hat{u}_k'}H_{\hat{f}_k'}\|_{L^2
(D_{2r})}^2<
e^{\|\hat{u}_k'\|_{L^\infty}}\int_{D_{2r}}|H_{\hat{f}_k'}|^2d\mu_{{\hat{f}_k'}}$,
we get $\|\nabla{\hat{f}_k'}\|_{W^{1,2}(D_{r})}<C$.
Recalling that $\hat{f}_k'$ converges in $C^0(D_r\setminus
D_\frac{r}{2})$, we complete the proof.

\endproof

\begin{rem}
In fact, we proved that $\mathcal{S}^*=\emptyset$.
\end{rem}

\section{Analysis of the neck}
For a sequence of conformal immersions from
a surface into $\R^n$ with the conformal class divergence,
the blowup comes from concentrations and collars.
Both cases can be changed into a blowup
analysis of a conformal immersion sequence
of $S^1\times[0,T_k]$ in $\R^n$ with $T_k\rightarrow+\infty$. So we first analyze the blow up procedure
on long cylinders without concentrations.

\subsection{Classification of bubbles of a simple
sequence over an infinite cylinder}

Let $f_k$ be an immersion sequence
of
$S^1\times [0,T_k]$ in $\R^n$ with $T_k\rightarrow+\infty$.
We say $f_k$ has concentration, if we can find a sequence
$\{(\theta_k,t_k)\}\subset S^1\times [0,T_k]$, such that
$$\lim_{r\rightarrow 0}\liminf_{k\rightarrow+\infty}\int_{D_r(\theta_k,t_k)}
|A_{f_k}|^2d\mu_{f_k}\geq 4\pi.$$
We say $\{f_k\}$ is simple if:
\begin{itemize}
\item[{\rm 1)}] $f_k$ has no concentration;
\item[{\rm 2)}] $f_k(S^1\times[0,T_k])$ can be extended to a compact closed
immersed surface $\Sigma_k$ with
$$\int_{\Sigma_k}(1+|A_{f_k}|^2)d\mu_{f_k}<\Lambda.$$
\end{itemize}

When $\{f_k\}$ is simple, we say $f_0$ is a bubble of $f_k$,
if we can find
a sequence $\{t_k\}\subset [0,T_k]$ with
$$t_k\rightarrow+\infty,\s and\s T_k-t_k\rightarrow+
\infty,$$
such that $f_0$ is a generalized
limit of $f_k(\theta,t_k+t)$. If $f_0$ is nontrivial,
we call it a nontrivial bubble.

For convenience, we call the generalized limit of $f(\theta,t+T_k)$
and $f(\theta,t)$ the top and the bottom
respectively.
Note that the top and the bottom  are in $W^{2,2}_{conf}(S^1\times(-\infty,
0])$ and $W^{2,2}_{conf}(S^1\times[0,+\infty))$
respectively.

\begin{defi}
Let $f^1$ and $f^2$ be two
bubbles which are limits of
$f_k(\theta,t+t_k^1)$
and $f_k(\theta,t+t_k^2)$
respectively.
We say these two bubbles are the same, if
$$\sup_k|t_k^1-t_k^2|<+\infty.$$
When $f^1$ and $f^2$ are not the same,
we say $f^1$ is in front of $f^2$ (or $f^2$ is behind $f^1$)
if  $t_k^1<t_k^2$. We say
$f^2$ follows $f^1$, if $f^2$ is behind $f^1$
and there are no non-trivial bubbles
between $f^1$ and $f^2$.
\end{defi}

Obviously,  the bubbles in this section must be
in $W^{2,2}_{conf}(S^1\times\R)$,
and must be one of the following:
\begin{itemize}
\item[1).] $S^2$-type, i.e. $I(f^0)(S^1\times\{\pm\infty\})
\neq 0$;
\item[2).] Catenoid-type, i.e. $I(f^0)(S^1\times\{\pm\infty\})=0$;
\item[3).] Plain-type, i.e. one and only one
of $I(f^0)(S^1\times\{\infty\})$, $I(f^0)(S^1\times\{-\infty\})$
is 0,
\end{itemize}
where $I=\frac{y-y_0}{|y-y_0|^2}$,
$y_0\notin f^0(S^1\times\R)$.

We give another classification of bubbles:
\begin{defi}\label{typeofbubble}
We call a  bubble
$f^0$ to be  a bubble of
\begin{itemize}
\item[] type $T_{\infty}$ if
$diam f^0(S^1\times\{+\infty\})=+\infty$; type $T_0$ if
$diam f^0(S^1\times\{+\infty\})=0$;
\item[] type $B_{\infty}$ if
$diam f^0(S^1\times\{-\infty\})=+\infty$; type $B_0$ if
$diam f^0(S^1\times\{-\infty\})=0$.
\end{itemize}
\end{defi}

We say $f_k$ has $m$ non-trivial bubbles, if we
can not find the ($m+1$)-th non-trivial bubble for any
subsequence of $f_k$.

\begin{rem}
Let  $f_0$ be a bubble. By \eqref{kappa2} and \eqref{kappa3},
$$\lim_{t\rightarrow+\infty}\int_{S^1\times\{t\}}\kappa_{f^0}
=2m^+\pi,\s and\s
\lim_{t\rightarrow+\infty}\int_{S^1\times\{t\}}\kappa_{f^0}
=2m^-\pi$$
for some $m^+$ and $m^-\in\mathbb{Z}$. Then
 $f^0$ is trivial implies that
 $\int_{S^1\times \R}K_{f^0}d\mu_{f^0}=0$. Thus
both $S^2$ type of bubbles and catenoid type of bubbles
are non-trivial.
\end{rem}

\begin{rem} It is easy to check that
$\mu(f^0)<+\infty$ implies that $f^0$ is
a sphere-type bubble and is of type  $(B_0,T_0)$.
\end{rem}

\begin{rem} If $f^{0'}$ is a catenoid-type bubble,
then it is of type $(B_\infty,T_\infty)$;
If $f^{0'}$ is a plain-type bubble,
then it is of type $(B_\infty,T_0)$
or  $(B_0,T_\infty)$.
\end{rem}

First, we study the case that
$f_k$ has no bubbles.
Basically, we want to show that after scaling, the
image of $f_k$ will converge to a topological disk.

\begin{lem}
If $f_k$ has no bubbles, then
$$\frac{diam\, f_k(S^1\times \{1\})}{diam\,
f_k(S^1\times\{T_k-1\})}\rightarrow 0\s
or\s +\infty.$$
\end{lem}

\proof Assume this lemma is not true. Then we may assume
$\frac{diam\, f_k(S^1\times \{1\})}{diam\,
f_k(S^1\times\{T_k-1\})}\rightarrow\lambda\in (0,+\infty)$. Let
$\lambda_k=diam f_k(S^1\times\{1\})$. By Theorem
\ref{convergence2},  $\frac{f_k(\theta,t)-f_k(0,1)}{\lambda_k}$
converges to $f^B$ weakly in $W^{2,2}_{loc}
(S^1\times(0,+\infty))$, and
$\frac{f_k(\theta,t+T_k)-f_k(0,T_k-1)}{\lambda_k}$ converges to
$f^T$ weakly in $W^{2,2}_{loc} (S^1\times(-\infty,0))$
respectively.

When $diam f^B(S^1\times \{+\infty\})=0$,  we set $\delta_k$ and
$t_k$ to be defined by
$$\delta_k=diam f_k(S^1\times\{t_k\})=
\inf_{t\in [1,T_k-1]}f_k(S^1\times\{t\}).$$
Obviously, $\delta_k\rightarrow 0$, and
 $t_k\rightarrow+\infty$, $T_k-t_k
\rightarrow+\infty$. $
\frac{f_k(\theta,t)-f_k(0,t_k)}{\delta_k}$
will converge to a non-trivial bubble. A contradiction.

When $diam f^B(S^1\times \{+\infty\})=+\infty$,
 we set $\delta_k'$ and $t_k'$ to be defined by
$$\delta_k'=diam f_k(S^1\times\{t_k'\})=
\sup_{t\in [1,T_k-1]}f_k(S^1\times\{t\}),$$
then we can also get a bubble.
\endproof

Now we assume $f_k$ has no
bubbles,
and $\frac{diam\, f_k(S^1\times \{1\})}{diam\,
f_k(S^1\times\{T_k-1\})}\rightarrow +\infty$.
Let $\lambda_k=diam\,
f_k(S^1\times\{T_k-1\})$. The bottom $f^B$  is the weak limit of
$f_k'=\frac{f_k(\theta,t)-f_k(0,1)}{\lambda_k}$.
Let $\phi$ be the conformal diffeomorphism
from $D\setminus\{0\}$ to $S^1\times[0,+\infty)$.
Then   $f^B\circ\phi$  is an immersion of $D$
in $\R^n$  perhaps with branch point $0$.
Moreover, by the arguments in \cite{C-L} or in \cite{C},
we have
$$f^B(\phi(0))=\lim_{t\rightarrow+\infty}\lim_{k\rightarrow+\infty}
f_k'(\theta,T_k-t).$$
Since $diam f_k'(S^1\times\{T_k-1\})\rightarrow 0$,
$f_k'(\theta,T_k-t)$  converges to a point, then the Hausdorff
limit of $f_k'((0,T_k))$ is a branched conformal immersion of $D$.

\begin{rem}
In fact, the above results and arguments hold for a sequence
 $\{f_k\}$ which has neither $S^2$-type
nor catenoid-type bubbles.\\
\end{rem}

Next, we show when $\{f_k\}$ has bubbles, how we will
find out all of them.
We need the following simple lemma:
\begin{lem}\label{interval} After passing to a subsequence, we can find
$0=d_k^0<d_k^1<\cdots<d_k^l=T_k$,
where $l\leq\frac{\Lambda}{4\pi}$,
such that
$$d_k^i-d_k^{i-1}\rightarrow+\infty,\s i=1,\cdots,l\s and
 \int_{S^1\times\{d_k^i\}}\kappa_k=2m_i\pi+\pi,\s m_i\in\mathbb{Z},\s
 i=1,\cdots,l-1,
$$
and
$$\lim_{T\rightarrow+\infty}\sup_{t\in [d_k^{i-1}+T,
d_k^i-T]}\left|
\int_{S^1\times \{t\}}\kappa_k-\int_{S^1\times\{d_k^{i-1}+T\}}\kappa_k\right|< \pi.$$
\end{lem}

\proof Let $\Lambda<4m\pi$. We prove the lemma by
induction of $m$.

We first prove it is true for $m=1$. Let
$$\lim_{t\rightarrow+\infty}
\lim_{k\rightarrow+\infty}\int_{S^1\times\{t\}}=2m_1\pi,\s
\lim_{t\rightarrow+\infty}
\lim_{k\rightarrow+\infty}\int_{S^1\times\{T_k-t\}}=2m_2\pi,$$
where $m_1$ and $m_2$ are integers.
Thus, we can find $T$, such that
$$\left|\int_{S^1\times\{T\}}\kappa_k-2m_1\pi\right|<\epsilon,\s
and \s \left|\int_{S^1\times\{T_k-T\}}\kappa_k-2m_2\pi\right|<\epsilon$$
when $k$ is sufficiently large. Take a $t_0\in (T,T_k-T)$,
such that
$$\int_{S^1\times[T,t_0]}|A_{f_k}|^2<2\pi,\s
\int_{S^1\times[t_0,T_k-T]}|A_{f_k}|^2\leq 2\pi.$$
By Gauss-Bonnet,
$$\left|\int_{S^1\times\{t\}}\kappa_k-\int_{S^1\times\{T\}}\kappa_k\right|
\leq \int_{S^1\times[T,t]}|K_{f_k}|d\mu_{f_k}
\leq \frac{1}{2}\int_{S^1\times[T,t_0]}|A_{f_k}|^2d\mu_{f_k}<\pi,
\s\forall t\in(T,t_0),$$
$$\left|\int_{S^1\times\{t\}}\kappa_k-\int_{S^1\times\{T_k-T\}}\kappa_k\right|
\leq
\frac{1}{2}\int_{S^1\times[t_0,T_k-T]}|A_{f_k}|^2d\mu_{f_k}<\pi,
\s\forall t\in(t_0,T_k-T).$$
Thus, we can take $\epsilon$ to be very small so that
$\int_{S^1\times\{t\}}
\neq 2i\pi$ for any $i\in\mathbb{Z}$ and
$t\in (T,T_k-T)$.

Now, we assume the result is true for $m$, and prove it
is also true for $m+1$. We have two cases.

Case 1, there is a sequence $\{t_k\}$, such that
$t_k\rightarrow+\infty$, $T_k-t_k\rightarrow+\infty$,
$\int_{S^1\times\{t_k\}}\kappa_k=2m_k\pi+\pi$
for some $m_k\in\mathbb{Z}$. For this case, we let
$f_k'=\frac{f_k(t+t_k,\theta)-f_k(t_k,0)}{\lambda_k}$
which converges weakly to $f_0'$, where $\lambda_k=diam f_k(S^1\times\{t_k\})$. Then by Gauss-Bonnet
$$\int_{S^1\times\R}|K_{f_0'}|\geq
\left|\int_{S^1\times(0,+\infty)}K_{f_0'}\right|
+\left|\int_{S^1\times(-\infty,0)}K_{f_0'}\right|\geq
2\pi.$$
Thus, $\int_{S^1\times\R}|A_{f_0'}|^2\geq 4\pi$. We can
find $T$, such that
$$\int_{S^1\times[0,t_k-T]}|A_{f_k}|^2<4(m-1)\pi,\s
and\s \int_{S^1\times[t_k+T,T_k]}|A_{f_k}|^2<4(m-1)\pi$$
when $k$ is sufficiently large. Thus, we can use
induction on $[0,t_k-T]$ to get $0=\bar{d}_k^0<
\bar{d}_k^1<\cdots<\bar{d}_k^{\bar{l}}=t_k-T$, and
on $[t_k+T,T_k]$ to get $t_k+T=\tilde{d}_k^0<\cdots<
\tilde{d}_k^{\tilde{l}}=T_k$. We can set
$$d_k^i=\left\{\begin{array}{ll}
            \bar{d}_k^i&i<\bar{l}\\
            t_k&i=\bar{l}\\
            \tilde{d}_k^{i-l}&i>\bar{l}
           \end{array}\right.$$
Then, we complete the proof.

\endproof

Set
$f_k^i=\frac{f_k(t+d_k^i,\theta)-f_k(d_k^i,0)}{
diam\, f_k(S^1\times\{d_k^i\})}$,
and assume $f_k^i\rightharpoonup f^i$.  It is easy to check that
$$\lim_{T\rightarrow+\infty}\lim_{k\rightarrow+\infty}
\int_{S^1\times\{d_k^i+T\}}\kappa_k
=\lim_{T\rightarrow+\infty}\lim_{k\rightarrow+\infty}
\int_{S^1\times\{d_k^{i+1}-T\}}\kappa_k,$$
we get
$$\lim_{T\rightarrow+\infty}\lim_{k\rightarrow+\infty}\int_{S^1\times[d_k^i+T,
d_k^{i+1}-T]}K_{f_k}=0.$$

\begin{rem}
In fact, we can get that for any $t_k<t_k'$ with
$$t_k-d_k^i\rightarrow+\infty, \s and\s d_k^{i+1}-t_k'\rightarrow+\infty,$$
we have
$$\lim_{k\rightarrow+\infty}\int_{S^1\times[t_k,
t_k']}K_{f_k}=0.$$
\end{rem}

Hence, we get

\begin{pro}\label{simple}
Let $f_k$ be a simple sequence on $S^1\times[0,T_k]$. Then after
passing to a subsequence,
$f_k$ has finitely many bubbles. Moreover, we have
$$
\lim_{T\rightarrow+\infty}\lim_{k\rightarrow+\infty}
\int_{S^1\times [T,T_k-T]}K_{f_k}d\mu_{f_k}=\sum_{i=1}^m
\int_{S^1\times\R}K_{f^i}d\mu_{f^i},$$
where $f^1$, $\cdots$, $f^m$ are all of the bubbles.
\end{pro}

Next, we prove a property of the order of the bubbles.

\begin{thm} Let $f^1$, $f^2$ be two bubbles. Then

1). If $f^1$ and $f^2$ are of type
$T_0$ and $B_0$  respectively, then
there is at least one catenoid-type bubble
between them.

2). If $f^1$ and $f^2$ are of type
$T_\infty$ and $B_\infty$  respectively, then
there is at least one $S^2$-type bubble
between $f^1$ and $f^{2}$.
\end{thm}

\proof 1).
Suppose $\frac{f_k(\theta,t_k^1+t)-f_k(0,t_k^1)}
{diam\, f_k(S^1\times\{t_k^1\})}\rightharpoonup f^1$, and
$\frac{f_k(\theta,t_k^{2}+t)-f_k(0,t_k^{2})}
{diam\, f_k(S^1\times\{t_k^2\})}\rightharpoonup f^{2}$.

Let $t_k'$  be defined by
\begin{equation}\label{infdiam}
diam\, f_k(S^1\times\{t_k'\})=\inf\{diam\,
f_k(S^1\times\{t\}):t\in[t_k^1+T,t_k^{2}-T]
\},
\end{equation}
where $T$ is sufficiently large. Since $f^1$ is
of type $T_0$ and  $f^{2}$  of type $B_0$,
we get
$$\lim_{t\rightarrow+\infty}diam\, f^1(S^1\times\{t\})=0,\s and\s
\lim_{t\rightarrow-\infty}diam\, f^{2}(S^1\times\{t\})=0.$$
Then, we have
$$t_k'-t_k^1\rightarrow+\infty,\s t_k^{2}-t_k'
\rightarrow+\infty.$$
If we set $f_k'(t)=\frac{f_k(\theta,t_k'+t)-
f_k(0,t_k')}{diam\, f_k(S^1\times\{t_k'\})}$,
then $f_k'$ will converge to a bubble $f'$
with
$$diam\, f'(S^1\times \{0\})
=\inf \{diam\,
f'(S^1\times\{t\}):t\in \R
\}=1.$$
Thus, $f'$ is a catenoid type bubble.

2). If we replace \eqref{infdiam} with
\begin{equation*}
diam\, f_k(S^1\times\{t_k'\})=\sup\{diam\,
f_k(S^1\times\{t\}):t\in[t_k^1+T,t_k^{2}-T]
\},
\end{equation*}
we will get 2).

\endproof

The structure of the bubble tree of a simple sequence
is  clear now: {\it The $S^2$ type bubbles  stand in
a line, with a unique catenoid type  bubble  between the two neighboring
$S^2$-type bubbles. There might
exist plain-type bubbles between the neighboring $S^2$ type
and catenoid type bubbles. A $T_0$ type bubble must follow a $B_\infty$ type bubble,
and a $T_\infty$ type bubble must follow a $B_0$ type bubble.}

\subsection{Bubble trees for a sequence of immersed $D$}

In this subsection, we will consider a conformal
immersion sequence
$f_k: D\rightarrow \R^n$ with $\mathcal{S}(f_k)=\{0\}$.
We assume that
$f_k(D)$ can be extended to a closed
embedded surface $\Sigma_k$
with
$$\int_{\Sigma_k}(1+|A_{\Sigma_k}|^2)d\mu<\Lambda.$$

Take $z_k$ and $r_k$, s.t.
\begin{equation}\label{top}
\int_{D_{r_k}(z_k)}|A_{f_k}|^2d\mu_{f_k}=4\pi-\epsilon,
\end{equation}
and $\int_{D_r(z)}|A_{f_k}|^2d\mu_{f_k}<4\pi-\epsilon$ for any $r<r_k$ and
$D_{r}(z)\subset D_\frac{1}{2}$, where $\epsilon$ is sufficiently small.

We set $f_k'=f_k(z_k+r_kz)-f_k(z_k)$. Then
$\mathcal{S}(f_k',D_L)=\emptyset$ for any $L$.
Thus, we can find $\lambda_k$, s.t.
$\frac{f_k'(z)}{\lambda_k}$ converges weakly to
$f^F$ which is a conformal immersion of $\C$
in $\R^n$. We call $f^F$ the first bubble of $f_k$
at the concentration point $0$.

It will be convenient to
make a conformal change of the domain. Let $(r,
\theta)$ be the polar coordinates centered at $z_k$.
Let $\varphi_k:S^1\times\R^1\rightarrow\R^2$ be the mapping
given by
$$r=e^{-t},\theta=\theta.$$
Then
$$\varphi_k^*(dx^1\otimes dx^1+dx^2\otimes dx^2)=
\frac{1}{r^2}(dt^2+d\theta^2).$$
Thus $f_k\circ\varphi_k$ can be considered as a conformal immersion
of $S^1\times [0,+\infty)$ in $\R^n$. For simplicity,
we will also denote $f_k\circ\varphi_k$ by $f_k$.

Set $T_k=-\log r_k$. Similarly to Lemma \ref{interval},
we have
\begin{lem}\label{interval2}
There is
$t_k^0=0<s_k^1<s_k^2< \cdots< s_k^l=T_k$,
such that $l\leq \frac{\Lambda}{4\pi}$ and

1). $\int_{S^1\times(s_k^i-1,s_k^i+1)}|A_{f_k}|^2\geq 4\pi$;

2). $\lim\limits_{T\rightarrow +\infty}\lim\limits_{
k\rightarrow+\infty}\sup\limits_{t\in [d_k^i+T,d_k^{i+1}-T]}
\int_{S^1\times(t-1,t+1)}|A_{f_k}|^2<4\pi$.

\end{lem}
Let
$f_k^i=f_k(\theta,s_k^i+t)$. A generalized limit of $f_k^i$
is called a bubble with concentration (which may be trivial).
There are $W^{2,2}$-conformal immersions
of $S^1\times\R$ with finite branch points and
finite $L^2$ norm of the second fundamental form.
However, if we neglect the concentration points,
we can also define the types of $T_\infty$, $T_0$, $B_{\infty}$,
and $B_0$ for it.

Obviously, we can find a $T'$, such that
$f_k$ is simple on  $S^1\times[s_k^{i}+T',s_k^{i+1}-T']$.
Note that the top  of $f_k$ on $S^1\times[s_k^i+T',s_k^{i+1}-T']$
is just a part of a generalized limit of $f_k^{i+1}$ and the
bottom of $f_k$ on $S^\times[s_k^i+T',s_k^{i+1}-T']$
is just a part of a generalized limit of $f_k^{i-1}$.
We call the union of nontrivial bubbles of $f_k$ on each $[s_k^i,s_k^{i+1}]$,
the generalized limit of $f_k^i$ and $f^F$
 the first level of bubble tree. By Proposition \ref{simple},
we have
$$\begin{array}{lll}
\lim\limits_{r\rightarrow 0}\lim\limits_{
k\rightarrow+\infty}\dint_{D_r}K_{f_k}&=&
\sum\limits_{i=1}^{l}\lim\limits_{T\rightarrow+\infty}
\lim\limits_{k\rightarrow+\infty}\dint_{S^1\times
[s_k^i-T'-T,s_k^i+T'+T]} K_{f_k^i}\\[\mv]
&&+
\sum\limits_{i=0}^l\lim\limits_{T\rightarrow+\infty}
\lim\limits_{k\rightarrow+\infty}\dint_{S^1\times
[s_k^i+T'+T,s_k^{i+1}-T'-T]} K_{f_k^i}\\[\mv]
&=&\sum\limits_{(r,\theta)\in \mathcal{S}(\{f_k^i\})}\lim\limits_{r\rightarrow 0}\lim\limits_{k\rightarrow+\infty}
\dint_{B_r(t,\theta)}K_{f_k^i}+
\sum_j\int_{S^1\times\R}K_{f^j},
\end{array}$$
where $\{f^j\}$ are all the bubbles of the first level.

Next,
at each concentration point of $\{f_k^i\}$, we get
the first level of $\{f_k^i\}$. We
usually call them the second level of bubble trees. Such a construction
will stop after finite steps.

\begin{lem}\label{identity1}After
passing to a subsequence,
$f_k$ has finitely many non-trivial bubbles.
Moreover, for any $r<1$
$$\lim_{k\rightarrow+\infty}
\int_{D_r}K_{f_k}d\mu_{f_k}=\int_{D_r}K_{f^0}d\mu_{f^0}+
\sum_{i=1}^m\int_{S^1\times\R}K_{f^i}d\mu_{f^i},$$
where $f^0$ is the generalized limit of $f_k$, and $f^1$, $f^2$,
$\cdots$, $f^m$ are all of the non-trivial
bubbles.
\end{lem}

\subsection{Immersion sequence of cylinder which is not simple}

Now we assume
$f_k$ is not simple on $S^1\times[0,T_k]$. We also assume
$f_k(S^1\times[0,T_k])$ can be extended to a closed
immersed surface $\Sigma_k$
with
$$\int_{\Sigma_k}(1+|A_{\Sigma_k}|^2)d\mu<\Lambda.$$
Moreover, we assume $f_k(t,\theta)$ and $f_k(T_k+t,\theta)$
have no concentration.

Then we still have Lemma \ref{interval2}.
The other properties
are the same as those of the immersion of $D$. Moreover, we have
$$\lim_{k\rightarrow+\infty}
\int_{S^1\times[0,T_k]}K_{f_k}d\mu_{f_k}=\int_{S^1\times[0,+\infty)}K_{f^B}
d\mu_{f^B}+
\int_{S^1\times(-\infty,0]}K_{f^T}d\mu_{f^T}+\sum_{i=1}^m\int_{\C}K_{f^i}d\mu_{f^i},$$
where $f^1$, $\cdots$, $f^m$
are all of the nontrivial bubbles.

\section{Proof of Theorem \ref{main}}
Since Theorem \ref{main2} can be deduced directly from
subsection 4.3, and Theorem \ref{sphere} can be deduced
directly from subsection 4.2,
 we only prove Theorem \ref{main}.

{\it Proof of Theorem \ref{main}:}
Take a curve $\gamma_i\subset\Sigma_0^i\setminus\mathcal{S}(\{f_k\circ
\psi_k\})$ with $\gamma_i(0)=p_i$. We set $\lambda_i=diam\,f_k(\gamma_i)$,
and $\tilde{f}_k^i=\frac{f_k\circ\psi_k-f_k\circ\psi_k(p_i)}{\lambda_k^i}$
which is a mapping from $\Sigma_0^i$
into $\R^n$. It is easy to
check that $\tilde{f}_k^i\in W^{2,2}_{{conf},loc}
(\Sigma_0^i,\psi_k^{\ast}(h_k),\R^n)$.

Given a point $p\in\Sigma_0^i$. We choose
 $U_k$, $U_0$, $\vartheta_k$, $\vartheta_0$ as in the
Theorem \ref{D.K.}.
Let $\hat{f}_k^i=\tilde{f}_k^i(\vartheta_k)$ which
is a map from $D$ into $\R^n$.
Let $V=\vartheta(D_\frac{1}{2})$. Since $\vartheta_k$
converges to $\vartheta_0$, $\vartheta_k^{-1}(V)\subset
D_\frac{3}{4}$
for any sufficiently large $k$.

When $p$ is not a concentration point, by Lemma \ref{measureconvergence},
for any
$\varphi$ with $supp\varphi\subset\subset V$,
we have
$$\int_{V}\varphi K_{\tilde{f}_k^i}d\mu_{\tilde{f}_k^i}
=\int_{D_\frac{3}{4}}\varphi(\vartheta_k)K_{\hat{f}_k^i}
d\mu_{\hat{f}_k^i}\rightarrow
\int_{D_\frac{3}{4}}\varphi(\vartheta_0)
K_{\hat{f}_0^i}=\int_V\varphi K_{f_0^i}d\mu_{f_0^i}.$$
When $p$ is a concentration point, by Lemma \ref{identity1}, we get
$$\int_{V}\varphi K_{\tilde{f}_k^i}d\mu_{\tilde{f}_k^i}
\rightarrow
\int_V\varphi K_{f_0^i}d\mu_{f_0^i}+
\varphi(p)\sum_j\int_{S^1\times \R}K_{f^i_j}d\mu_{f^i_j},$$
where  $\{f^i_j\}$ is the set of  nontrivial bubbles of $\hat{f}_k^i$
at $p$.

Next, we consider the convergence of $f_k$ at
the collars. Let $a^j$ be the intersection of $\overline{\Sigma_0^i}$
and $\overline{\Sigma_0^{i'}}$.
We set $\check{f}_k^j=f_k(\phi_k^j)$, and $T_k^j
=\frac{\pi^2}{l_k^j}-T$.
We may choose $T$ to be sufficiently large such that
$\check{f}_k^j(T_k^j-t,\theta)$ and $\check{f}_k^j
(-T_k^j+t,\theta)$ have no blowup point. Then $\check{f}_k^j$
satisfies the conditions in subsection 2.4.
So the convergence of $\check{f}_k^j$ is clear.
Since
$$\check{f}_k^j=f_k\circ\phi_k^j=f_k\circ\psi_k\circ(\varphi_k\circ\phi_k^j)=
\tilde{f}_k(\varphi_k\circ\phi_k^j).$$
The images of the limit of $\check{f}_k^j(T_k^j-t,\theta)$ and $\check{f}_k^j
(-T_k^j+t,\theta)$
are parts of the images of $\tilde{f}_0^i$ and $\tilde{f}_0^{i'}$.
Then we have
$$\lim_{\delta\rightarrow 0}
\lim_{k\rightarrow+\infty}\int_{\Sigma_0(\delta,a^j)}
K_{f_k}=\sum_i\int_{S^1\times\R}K_{f^{i'}},$$
where all $f^{i'}$ are  nontrivial bubbles of $\check{f}_k^j$.

\endproof

\section{A remark about trivial bubbles}
The methods in section 4 can be also used to find all bubbles
with $\|A\|_{L^2}\geq\epsilon_0$ for a fixed $\epsilon_0>0$.
We only consider the simple sequence $f_k$ on
$S^1\times[0,T_k]$ here.

Let $t_k$ be a sequence with $t_k, T_k-t_k\rightarrow \infty$, such that
$\frac{f_k(t+t_k,\theta)-f_k(t_k,0)}{\lambda_k}$
converges to a $f_0\in W^{2,2}(S^1\times\R,\R^n)$ with
$\int_{S^1\times\R}|A_{f_0}|^2\geq\epsilon_0^2$.
Take $T$, such that $\int_{S^1\times[-T,T]}
|A_{f_0}|^2\geq\frac{\epsilon_0^2}{2}$. We consider
 the convergence on $S^1\times [0,t_k-T]$
  and $S^1\times[t_k+T,T_k]$ respectively.
In this way, we can find out all the bubbles.

\end{document}